\documentclass{article}
\usepackage{times}
\usepackage{amsmath}
\usepackage{amssymb}
\usepackage{amsthm}
\usepackage{graphicx}
\usepackage{color}
\usepackage{wrapfig}

\newcommand {\mm}[1] {\ifmmode{#1}\else{\mbox{\(#1\)}}\fi}
\renewcommand{\proof}{\noindent{\sc Proof.~}}
\newcommand{\eop}{\hfill\usebox{\smallProofsym}\bigskip}  %
\newsavebox{\smallProofsym}                            
\savebox{\smallProofsym}                               %
{
\begin{picture}(6,6)                                   %
\put(0,0){\framebox(6,6){}}                            %
\put(0,2){\framebox(4,4){}}                            %
\end{picture}                                          %
}

\makeatletter
\long\def\@makecaption#1#2{%
  \vskip\abovecaptionskip
  \sbox\@tempboxa{\small #1: #2}%
  \ifdim \wd\@tempboxa >\hsize
    \small #1: #2\par
  \else
    \global \@minipagefalse
    \hb@xt@\hsize{\hfil\box\@tempboxa\hfil}%
  \fi
  \vskip\belowcaptionskip}
\makeatother

\newcommand{\Nspace}        {\mm{{\mathbb N}}}
\newcommand{\Rspace}        {\mm{{\mathbb R}}}
\newcommand{\Sspace}        {\mm{{\mathbb S}}}

\newcommand{\Gnet}[1]       {\mm{{\mathcal G}_{#1}}}
\newcommand{\Hnet}[1]       {\mm{{\mathcal H}_{#1}}}
\newcommand{\cell}[1]       {\mm{\rm cell}{({#1})}}
\newcommand{\ucell}[1]      {\mm{\rm ucell}{({#1})}}

\newcommand{\dist}[2]       {\mm{\|{#1}-{#2}\|}}
\newcommand{\Length}[1]     {\mm{\ell}{({#1})}}
\newcommand{\scalprod}[2]   {{\langle #1 , #2 \rangle}}

\newcommand{\mmm}           {\mm{\bf m}}
\newcommand{\ddd}[1]        {\mm{\bf d}_{#1}}

\newtheorem{definition}{{\sc Definition}}
\newtheorem{theorem}{{\sc Theorem}}
\newtheorem{lemma}{{\sc Lemma}}
\newtheorem{proposition}{{\sc Proposition}}

\newcommand{\ignore}[1]  {}

\begin{document}

\begin{center}
  {\bf \Large On the Configuration Space of Steiner Minimal Trees}
    \footnote{The authors acknowledge the support by the Russian Government
              through the mega project, resolution no.\ 220,
              contract no.\ 11.G34.31.0053.}

        H.~Edelsbrunner\footnote{IST Austria (Institute of Science
          and Technology Austria), Klosterneuburg, Austria.} and
        N.\,P.~Strelkova\footnote{Moscow State University,
          Moscow, Russian Federation.}
\end{center}

\begin{abstract}
  Among other results, we prove the following theorem about
  Steiner minimal trees in $d$-dimensional Euclidean space:
  if two finite sets in $\Rspace^d$ have unique and combinatorially
  equivalent Steiner minimal trees, then there is a homotopy between
  the two sets that maintains the uniqueness and the combinatorial
  structure of the Steiner minimal tree throughout the homotopy.
  

05C05, 51M16 
\end{abstract}

\section{Introduction}
\label{sec1}

Shortest networks connecting finite sets of geographic locations were
intensely studied when the telephone providers were legally required
to charge for phone calls an amount that is proportional to the
length of the connection \cite{CLDN04}.
This motivated the discovery of many mathematical and computational
properties of such networks.
Nevertheless, many of the fundamental questions are still open.
For example, in 1968 Gilbert and Pollak \cite{GiPo68} conjectured
that the ratio of the length of a minimum spanning tree over
the length of a shortest network of the same point set
cannot exceed $2 / \sqrt{3}$.
In 1992, Du and Hwang claimed a proof of the conjecture \cite{DuHw92},
but $10$ years later, Ivanov and Tuzhilin noted that there
are gaps \cite{IvTu03}, and the problem remains open until today \cite{IvTu12}.
In contrast to the ratio problem, surprisingly little attention was
directed toward the configuration space of the shortest networks,
a topic we approach in this paper.
Specifically, we consider a finite set of points, $S$, in $d$-dimensional
Euclidean space, which we denote as $\Rspace^d$.
A \emph{spanning network} of $S$ is a finite set of rectifyable curves,
each connecting two points in $\Rspace^d$ but not necessarily in $S$,
whose union is connected and contains $S$.
We measure the length of the network using the Euclidean metric.
A \emph{shortest network} of $S$ is a spanning network of minimum length.
\begin{wrapfigure}[10]{r}{120pt}
  \centering \resizebox{!}{1.3in}{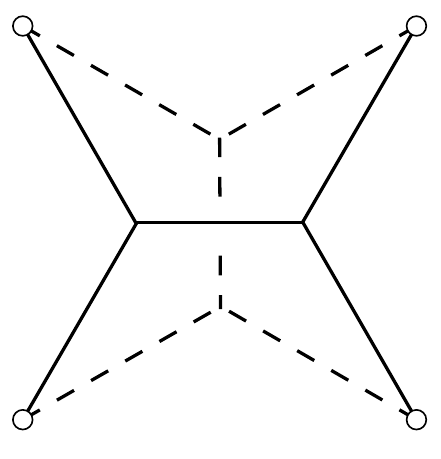}
  \caption{Two Steiner minimal trees of the four corners of a square.}
  \label{fig:square}
\end{wrapfigure}
It is known to exist, and it satisfies the following properties;
see e.g.\ \cite{GiPo68}:
\begin{itemize}
  \item it is a tree whose edges are straight line segments ending
    at vertices of degree $1$, $2$, or $3$;
  \item all degree-$1$ and degree-$2$ vertices are points in $S$;
  \item the angles between the edges meeting at a degree-$3$ vertex
    are $120^{\circ}$, while the angle between the edges meeting
    at a degree-$2$ vertex is at least $120^{\circ}$.
\end{itemize}
A shortest network of a finite set of points is usually referred
to as a \emph{Steiner minimal tree}.
This tree may or may not be unique, and in the latter case,
we call the Steiner minimal tree \emph{ambiguous}.
For example, the four corners of a square have two Steiner minimal trees;
see Figure \ref{fig:square}.
Note that these two trees have different combinatorial structure if the
points are labeled.
Just to mention one difference: the points $A$ and $B$ are connected
by a path of two edges in the solid tree, and by a path of three edges
in the dashed tree.
Suppose now that $S$ consists of four (arbitrary) points in the plane,
and we need to find a Steiner minimal tree.
If we know how many additional vertices there are in the tree,
and how they are connected to the four given points and to each other,
then we can locate them using the $120^{\circ}$-condition.
Generalizing this idea to $n$ points in the plane gives Melzak's
algorithm \cite{Mel61}.
But a priori, we do not know which combinatorial structure gives
the minimum length,
and there is no known way that avoids checking all of them in the worst case.
Indeed, the problem of constructing a Steiner minimal tree
in $\Rspace^d$ is NP-hard even for $d=2$; see \cite[page 209]{GaJo79}.

For a given finite set, $S$, the collection of Steiner minimal trees
of $S$ is well defined, and every tree in this collection has a different
\emph{combinatorial type} \cite{GiPo68},
a notion that will be made precise later.
We are interested in the following question:
{\bf for a given combinatorial type, what does the space of points
     whose Steiner minimal trees are of this type look like topologically?}
To make this more concrete, we order the $n$ points in $S \subseteq \Rspace^d$
and list the coordinates in sequence,
mapping $S$ to a point in $\Rspace^{nd}$.
For a given combinatorial type, we consider the space of points
in $\Rspace^{nd}$ that correspond to sets
in $\Rspace^d$ with unambiguous Steiner minimal tree of the given type.
Our main result is that this space is path-connected.
We have a weaker result if we relax the requirement and allow
ambiguous Steiner minimal trees.
Here we can prove path-connectivity only for full networks in
the plane.
Finally, we exhibit an example to show that the space of points
with necessarily ambiguous Steiner minimal trees of a given type
is generally not connected.
While the restricted setting for the result on cells is
disappointing, it is consistent with a general lack of understanding
of ambiguous Steiner minimal trees.
For example, it is not known that the space of ambiguous
Steiner minimal trees has measure zero.
The only partial result is in $d = 2$ dimensions, where 
the unambiguous Steiner minimal trees
contain an everywhere dense open subset of the configuration space,
but the proof in \cite{IvTu06} is long and complicated.
The approach in \cite{Obl09} leads to a shorter proof.

\paragraph{Outline.}
Section \ref{sec2} introduces the main concepts needed to formally
state our results.
Section \ref{sec3} presents the key tools used in our proofs.
Section \ref{sec4} presents the proofs.
Section \ref{sec5} concludes the paper.

\section{Concepts and Results}
\label{sec2}

The two main concepts in this paper are the spanning networks of
a set of $n$ points in $\Rspace^d$,
and the decomposition of the $nd$-dimensional configuration space
into cells whose points have combinatorially equivalent shortest networks.
We discuss both and introduce the terminology needed to formally
state our results.

\paragraph{Graphs, networks, types.}
Similar to the notion of a curve in differential geometry,
we define a network as a map from an abstract topological
to a concrete geometric space.
In our context, the abstract space is a finite connected graph,
$G = (V \sqcup W, E)$.
It is \emph{partially ordered}, with $V$ ordered and $W$ unordered.
The difference arises because the vertices in $V$ are mapped
to the points in $S$, which are ordered,
while the vertices in $W$ can be mapped anywhere.
To emphasize the difference, we sometimes call an element of $V$
a \emph{terminal vertex} and an element of $W$ an \emph{interior vertex};
an element of $E$ is an \emph{edge}.
Two partially ordered graphs,
$G = (V \sqcup W, E)$ and $G' = (V' \sqcup W', E')$,
are \emph{combinatorially equivalent} if there exist compatible
bijections $V \to V'$, $W \to W'$, and $E \to E'$
of which the first preserves order.
Observe that the preservation of order prevents the existence
of such bijections for the two trees in Figure \ref{fig:square}.
We have $V = V' = \{A, B, C, D\}$,
and since we are obliged to map $A$ to $A$ and so on,
it is impossible to find a compatible bijection between the two
edge sets.

We find it convenient to consider a graph as a $1$-dimensional
topological space, so it makes sense to talk about a point $x \in G$,
which may be a vertex or a point on an edge of $G$.
If $G$ and $G'$ are combinatorially equivalent, then there exists
a \emph{consistent} homeomorphism $g \colon G \to G'$,
that is: $g$ restricts to the bijections
establishing the combinatorial equivalence.
With this preamble, we introduce the first main concept.
\begin{definition}
  A \emph{network} is a continuous map $\Gnet{} \colon G \to \Rspace^d$
  such that the image of every edge is a rectifyable curve or a point,
  and the preimage of the image of every vertex is a subgraph of $G$,
  possibly together with subsets of points on edges
  whose images happen to pass through this image.
  We call $G$ a \emph{parametrization} of $\Gnet{}$,
  and a \emph{contracting parametrization} if every subgraph with constant
  image is connected.
\end{definition}
The \emph{(combinatorial) type} of $\Gnet{}$  is the graph $G_0$
obtained from $G$ by shrinking each subgraph with constant image to a vertex.
Shrinking is similar to the contraction operation defined
in Tutte \cite[page 32]{Tut01} but more general because we do not
require that the subgraphs that shrink to vertices be connected.
The type is again a graph although it may have loops and multi-edges
even if $G$ does not.
The network $\Gnet{} \colon G \to \Rspace^d$ \emph{spans}
an ordered set $S \subseteq \Rspace^d$ if $\Gnet{}$
restricts to an order-preserving bijection from $V$ to $S$.
To compare networks, we 
write $\Length{\Gnet{}}$ for the sum of the lengths of the curves.
A \emph{straight-line network} maps every edge to a line segment or point.
Of particular importance are shortest spanning networks,
which are necessarily straight-line but also satisfy the three
properties stated in Section \ref{sec1}.
A \emph{Steiner minimal tree} is the image of a shortest network.
The image of a shortest network is necessarily a tree,
so it is natural to parametrize it using a tree.
In this case, the type is obtained by shrinking connected subtrees,
so the parametrization is necessarily contracting.

The above definitions allow for different degrees of sameness between
two networks,
$\Gnet{} \colon G \to \Rspace^d$ and $\Gnet{}' \colon G' \to \Rspace^d$.
We call $\Gnet{}, \Gnet{}'$ \emph{identical} if $G, G'$ are
combinatorially equivalent and there is a consistent homeomorphism
$g \colon G \to G'$ such that $\Gnet{} (x) = \Gnet{}' (g(x))$
for every $x \in G$.
Note that this is strictly stronger than requiring combinatorially
equivalent graphs and equal images.
A weaker and perhaps more natural notion calls
$\Gnet{}, \Gnet{}'$
\emph{the same} if their types, $G_0, G_0'$, are combinatorially equivalent,
and there is a consistent homeomorphism $g_0 \colon G_0 \to G_0'$
such that $\Gnet{} (x) = \Gnet{}' (g_0(x))$ for every $x \in G_0$.
For example, two shortest networks are the same iff
they have the same Steiner minimal tree.
Yet weaker is the notion that says $\Gnet{}, \Gnet{}'$
\emph{have the same type} if $G_0, G_0'$ are combinatorially equivalent.
In dimension $d = 2$, we require in addition that there is a
homeomorphism $h \colon \Rspace^2 \to \Rspace^2$ such that
$\Gnet{} (x) = h (\Gnet{}' (g_0(x)))$.
For plane embeddings, this is equivalent to requiring that
the ordering of the edges around the vertices is the same in both images.
While the three notions of sameness are different and relate
to each other by a chain of implications, they are all transitive.
We get the reverse chain of implications for the negations:
having different type implies not being the same implies not being identical.

\paragraph{Configurations, cells, decompositions.}
Let $S$ be an ordered set of points
$p_i = (p_{i1}, p_{i2}, \ldots, p_{id}) \in \Rspace^d$ for $1 \leq i \leq n$.
Listing the coordinates in sequence, we obtain a point
\begin{align}
  s  &=  (p_{11}, p_{12}, \ldots, p_{1d},
          p_{21}, p_{22}, \ldots, p_{nd} ) \in \Rspace^{nd} .
\end{align}
We require that the points in $S$ are distinct,
which implies that $s$ cannot lie on any of the
$(nd - d)$-dimensional planes of the form
$\Delta_i^j = \{ (\ldots, x, \ldots, x, \ldots) \in \Rspace^{nd}
              \mid x \in \Rspace^d \}$,
where we list the $i$-th and the $j$-th $d$-tuple of coordinates,
both specifying the same point in $\Rspace^d$.
Observe that all $\Delta_i^j$ contain the diagonal of points $(x,x,\ldots,x)$.
We use $\Rspace^{nd}_\Delta = \Rspace^{nd} \setminus \bigcup_{i < j} \Delta_i^j$
as the \emph{configuration space} for sets of $n$ points in $\Rspace^d$.
\begin{definition}
  The \emph{cell} of a type $G$, denoted as $\cell{G}$,
  is the space of points $s \in \Rspace^{nd}_\Delta$
  that have shortest networks of type $G$.
  The \emph{unambiguous cell} of $G$, denoted as $\ucell{G}$,
  is the subspace of points $s \in \cell{G}$ such that $s$
  has a unique Steiner minimal tree.
\end{definition}
The remaining points, in $\cell{G} \setminus \ucell{G}$,
have at least two shortest networks,
one of type $G$ and another of necessarily different type;
see \cite{GiPo68}.
What can we say about the cells?
Clearly, they cover the configuration space
since every point $s \in \Rspace^{nd}_\Delta$ has at least
one shortest network.
The cells have possibly non-empty intersections,
and removing these intersections leaves us with the unambiguous cells,
which are disjoint but do not necessarily cover the configuration space.
As mentioned earlier, there is not much known about the space of
points with ambiguous shortest networks.
To complete the picture, we note that the points in the $\Delta_i^j$
correspond to multisets, or after identification to
sets of size $m < n$ in $\Rspace^d$.
We can therefore complete the decomposition of $\Rspace^{nd}$
by partitioning the $\Delta_i^j$ into strata within which $m$ is a constant,
and to decompose each stratum according to shortest networks.
By definition, the cells in these strata belong to the
boundaries of the cells in $\Rspace^{nd}_\Delta$,
but not to the cells themselves.

\begin{figure}[hbt]
  \centering \resizebox{!}{1.2in}{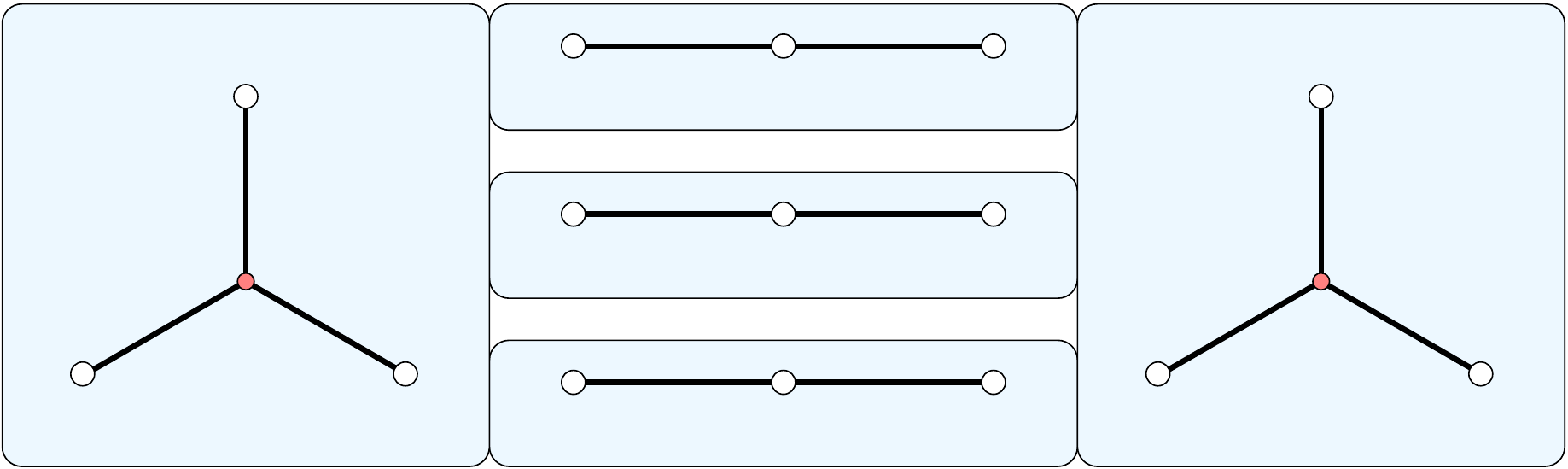}
  \caption{The five types of shortest networks for three points
    in the plane.
    Two boxes touch iff there is a transition between the
    types that avoids configurations with fewer than three points.}
  \label{fig:example}
\end{figure}
We illustrate the decomposition of the configuration space
with a small example: $n = 3$ points in $d = 2$ dimensions.
There are only five different types, each generating a cell
in $\Rspace^6_\Delta$; see Figure \ref{fig:example}.
Factoring out the three dimensions of rigid motions,
we may visualize the decomposition in the remaining $3$ dimensions.
Drawing the diagonal as a point at the origin
with three emanating half-lines for the
$\Delta_i^j$, $1 \leq i < j \leq 3$,
we get a ring of three cells separating the two remaining cells.
The origin lies in the boundaries of all five cells,
while each half-line belongs to the boundaries of four.
Only six of the ten transitions between the five cells avoid
the origin and the half-lines,
and we visualize them as $2$-dimensional membranes,
which lie in the common boundaries of the corresponding cells.
In this particular example, these membranes belong to the three cells
in the ring but not to the other two cells.
For example, a point on the membrane separating the cell on the left
from the cell at the center in Figure \ref{fig:example}
has a Steiner minimal tree consisting of two edges
that form an angle of $120^{\circ}$ at $B$.
The combinatorial type of this network is different from that on the left
but the same as that at the center.
Indeed, the two separated cells do not contain any points of their boundary,
while the three cells in the ring contain every boundary point
that does not belong to the origin and the half-lines.

\paragraph{Main results.}
We are now ready to give formal statements of our main results.
They state that unambiguous cells are connected for all dimensions $d$,
and that cells are connected for full networks in $\Rspace^2$.
In contrast, the difference between a cell and the corresponding
unambiguous cell is not necessarily connected,
which we will show by exhibiting an example of 
$n = 4$ points in $d = 2$ dimensions.
\begin{theorem}\label{thm:connected-ucells}
  {\sc (Connectedness of unambiguous cells)}
  Let $G$ be the type of a shortest network of a finite set of
  points in $\Rspace^d$, for $d \geq 2$.
  Then $\ucell{G}$ is path-connected.
\end{theorem}
\begin{theorem}\label{thm:connected-cells}
  {\sc (Connectedness of cells)}
  Let $G$ be the type of a shortest network of a finite set of
  points in $\Rspace^2$
  such that all terminal vertices have degree-$1$.
  Then $\cell{G}$ is path-connected.
\end{theorem}

\section{Tools}
\label{sec3}

\begin{wrapfigure}[14]{r}{2.3in}\vspace*{-1.3cm}
  \centering \resizebox{!}{1.8in}{\includegraphics{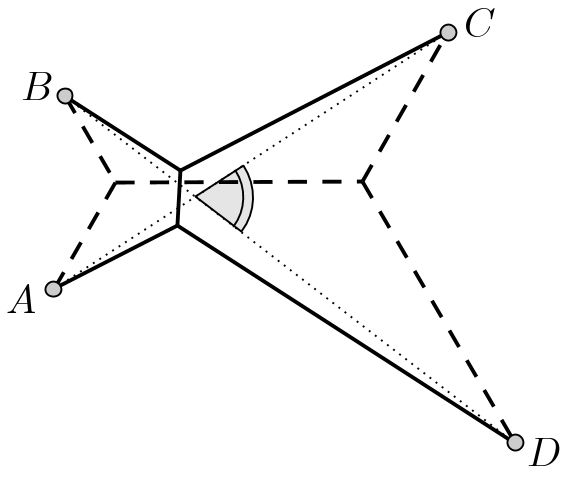}}
  \caption{Two locally minimal networks spanning the same four points
    in the plane.
    The dashed tree is shorter because its center edge aligns with
    the smaller angle pair formed by the diagonals;
    see \cite[Theorem 6]{Pol78}.}
  \label{fig:smt4points}
 \end{wrapfigure}

The proof of our main result is inductive, applying a procedure
similar to Melzak's algorithm;
see \cite[Section 5]{GiPo68} or \cite{Mel61}.
We begin by introducing the ingredients of the procedure,
proving lemmas in preparation of the proofs of the theorems.

\paragraph{Locally minimal and pairs of codirected networks.}
A network that satisfies the three properties stated in Section \ref{sec1}
is called \emph{locally minimal}.
Every shortest network is locally minimal
but there are others; see Figure \ref{fig:smt4points} for an example.
As its name suggests, there are no shorter networks nearby,
in which the notion of nearness is defined by parametrizing the
networks by their interior vertices.
More generally, it is not possible to have two locally minimal
networks of the same type.
We state this result for later reference;
see \cite[pages 927 and 939]{IvTu99}.
\begin{proposition}\label{prop:minimal-networks}
  {\sc (Uniqueness of minimal network)}
  Let $G$ be a partially ordered tree with $n$ terminal vertices
  and $S$ an ordered set of $n$ points in $\Rspace^d$.
  If there exists a locally minimal network that spans $S$ and is
  parametrized by $G$, then it is unique and the shortest such network.
\end{proposition}
As a special case, we see that two different shortest networks
of the same set $S$ that are parametrized by the same tree are impossible.
A related concept is the following.
Two straight-line networks, $\Gnet{} \colon G \to \Rspace^d$ and
$\Hnet{} \colon H \to \Rspace^d$, that span the same set $S$
are \emph{codirected} if they look geometrically the same in the
neighborhood of every vertex.
More formally, we require that there exists $r > 0$ such that the
intersection of the ball with radius $r$ centered at
$\Gnet{} (u) = \Hnet{} (u)$ with the images of $\Gnet{}$ and $\Hnet{}$
are the same for every terminal vertex $u$.
The existence of codirected locally minimal trees that are not the same is open,
but in the plane it is known that such pairs do not exist \cite{Obl09}.
\begin{proposition}\label{prop:codirected}
  {\sc (Non-existence of codirected trees)}
  Any two locally minimal trees spanning the same ordered set
  of finitely many points in $\Rspace^2$ that are codirected are the same.
\end{proposition}

\paragraph{Sequences of networks.}
Another useful tool in the study of cells and their boundaries
are converging sequences.
We write $\Gnet{i}$ for the networks in an infinite sequence,
indexing them with the positive integers $i \in \Nspace$.
In this paper, we are primarily interested in sequences such that the
$\Gnet{i}$ have all the same type, $G$.
We say such a sequence $\Gnet{i}$ \emph{converges} if the sequence of points
$\Gnet{i} (u)$ converges for every vertex $u$ of $G$.
For straight edges, the limit of a converging sequence is
well defined, namely the straight-line network $\Gnet{} \colon G \to \Rspace^d$
such that $\Gnet{i} (u) \to \Gnet{} (u)$ for each vertex $u$.
We write $\Gnet{i} \to \Gnet{}$ as well as
$\Gnet{} = \lim_{i \to \infty} \Gnet{i}$.
Note that $\Gnet{}$ is not necessarily of type $G$,
but $G$ parametrizes $\Gnet{}$ by construction.
An important method for obtaining converging sequences
follows from the Bolzano-Weierstrass Theorem in analysis:
starting with an infinite sequence of networks with finitely many
types and points within a compact domain,
we obtain a converging sequence by repeatedly restricting
the choice and continuing the process with an infinite subsequence.

For straight-line networks,
$\Gnet{i} \to \Gnet{}$ implies that the
length of $\Gnet{i}$ converges to the length of $\Gnet{}$.
It is therefore easy to see that if all $\Gnet{i}$ are shortest networks,
then $\Gnet{}$ is a shortest network with a contracting parametrization.
Another important property of converging sequences of shortest
networks is that convergence for the terminal vertices
implies convergence for all vertices.
In other words, as long as the shortest network does not
change its type, it depends continuously on the points it spans.
\begin{lemma}\label{lemma:convergence}
  {\sc (Convergence of shortest networks)}
  Let $G$ be a tree, and let $s_i$ be a converging sequence of points
  in $\cell{G}$ with shortest networks $\Gnet{i} \colon G \to \Rspace^d$,
  all of type $G$.
  Then $\Gnet{i}$ is a converging sequence of networks,
  and the limit network is a shortest network of the limit of $s_i$.
\end{lemma}
\proof
 Since $s_i$ converges, there exists a bounded and convex domain
 $D \subseteq \Rspace^d$ that contains the images of the terminal
 vertices of $G$ for all $i \in \Nspace$.
 The images of the interior vertices lie within the convex hull
 of the images of the terminal vertices and hence also in $D$.
 We can therefore apply the Bolzano-Weierstrass Theorem
 and obtain a converging subsequence of shortest networks,
 with limit $\Gnet{}$.
 To get a contradiction to the claimed property,
 we assume that $\Gnet{i}$ does not converge to $\Gnet{}$.
 Then there exists an interior vertex $u$ of $G$, a constant $C > 0$,
 and an infinite subsequence $\Gnet{i_j}$
 such that $\dist{\Gnet{i_j} (u)}{\Gnet{} (u)} > C$ for all $j$.
 Applying the Bolzano-Weierstrass Theorem to this subsequence
 gives another converging sequence of shortest networks,
 now with limit $\Gnet{}' \neq \Gnet{}$.
 Thus, we get two shortest networks of the same set $S \subseteq \Rspace^d$,
 both parametrized by the same tree $G$.
 Both are locally minimal, which contradicts Lemma \ref{prop:minimal-networks}.
\eop

Note that Lemma \ref{lemma:convergence} implies that every point $s$
in the boundary of $\cell{G}$ has a shortest network parametrized by $G$,
and the parametrization is contracting.

\paragraph{Moustaches.}
A particularly simple application of infinite sequences is used to shrink
and grow edges of Steiner minimal trees that end at degree-$1$ vertices.
Let $\Gnet{}$ be a shortest network of a set $S$ of $n \geq 2$
points in $\Rspace^d$, let $G$ be its type,
and recall that it enjoys the three properties stated in Section \ref{sec1}.
Hence, $G$ has at least one of the following two substructures:
a \emph{two-sided moustache} consisting of degree-$1$ vertices $p$ and $q$
both adjacent to a degree-$3$ vertex $v$,
or a \emph{one-sided moustache} consisting of a degree-$1$ vertex $p$
adjacent to a degree-$2$ vertex $v$.
We call $v$ the \emph{anchor} of the moustache in both cases.
We \emph{shave} a two-sided moustache by removing $p$ and $q$ together
with the edges connecting them to $v$.
This leaves a smaller graph $G'$ and a network
$\Gnet{}' \colon G' \to \Rspace^d$ that spans a smaller set of points,
namely $S' = S \setminus \{ \Gnet{} (p), \Gnet{} (q) \}$
if $\Gnet{} (v) \in S$, and union $\{ \Gnet{} (v) \}$ if $\Gnet{} (v) \not\in S$.
Similarly, we shave a one-sided moustache by removing $p$ together
with the edge connecting it to $v$.
In this case, $\Gnet{} (v)$ is necessarily in $S$,
so we get a spanning network of $S' = S \setminus \{ \Gnet{} (p) \}$.
Importantly, the operation preserves minimality in both cases.
For later reference, we state and prove a slightly more general claim,
namely that trimming a moustache preserves minimality.
\begin{lemma}\label{lemma:trimming}
  {\sc (Trimming a moustache)}
  Let $G$ be the type of a shortest network $\Gnet{}$ of a set of $n \geq 2$
  points in $\Rspace^d$, let $p$ be a degree-$1$ vertex of $G$
  and $v$ its neighbor.
  For $0 \leq t \leq 1$,
  let $\Gnet{t} \colon G \to \Rspace^d$ be the same except that it
  maps $p$ to $\Gnet{t} (p) = (1-t) \Gnet{} (p) + t \Gnet{} (v)$.
  Then $\Gnet{t}$ is a shortest network of
  $S_t = S \setminus \{ \Gnet{} (p) \} \cup \{ \Gnet{t} (p) \}$,
  and if $\Gnet{}$ is unambiguous then so is $\Gnet{t}$.
\end{lemma}
\proof
 Suppose there exists a network $\Hnet{t} \colon G \to \Rspace^d$
 that spans $S_t$ with $\Length{\Hnet{t}} < \Length{\Gnet{t}}$.
 Append the line segment from $\Hnet{t} (p) = \Gnet{t} (p)$
 to $\Gnet{} (p)$ to get a network $\Hnet{}$ that now spans $S$.
 Then $\Length{\Hnet{}} < \Length{\Gnet{}}$,
 which contradicts the assumption that $\Gnet{}$ is a shortest
 network of $S$.
 The same argument contradicts the existence of different shortest
 networks for $S_t$ if there is only one for $S$.
\eop

\begin{figure}[hbt]
  \centering \resizebox{!}{1.6in}{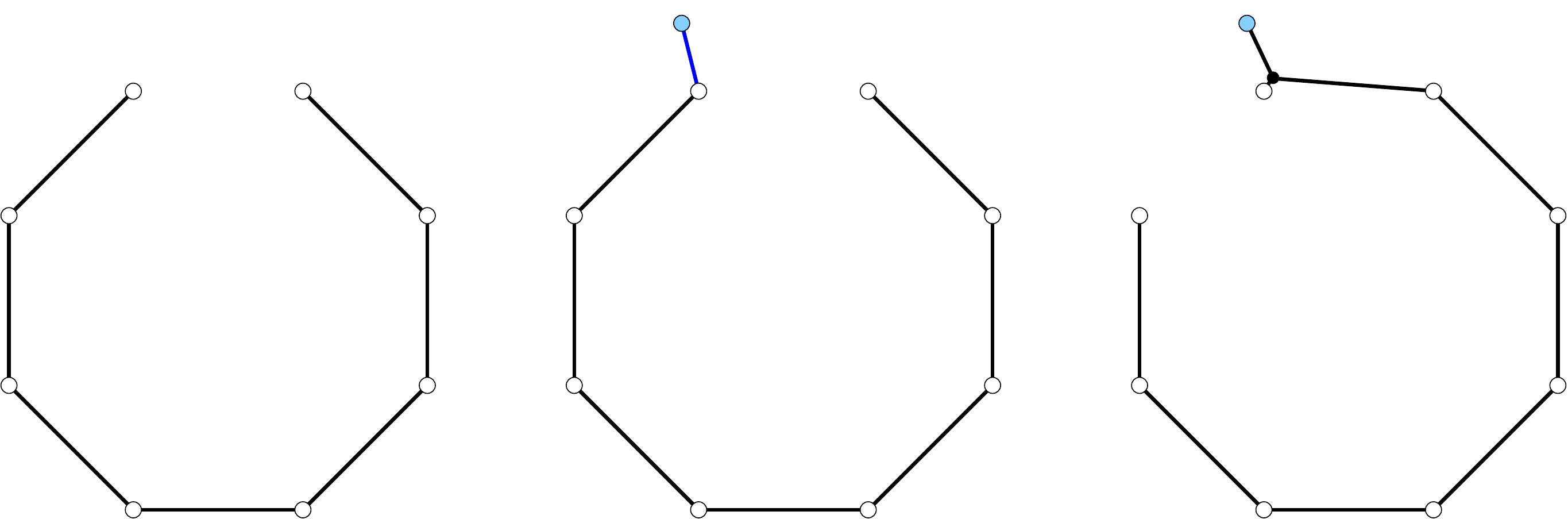}
  \caption{From \emph{left} to \emph{right}:
    a Steiner minimal tree of the corners of a regular octagon,
    the tree with a one-sided moustache attached,
    a shorter network with a degree-$3$ vertex connecting the
    endpoint of the moustache to the rest of the tree.}
  \label{fig:octagon}
\end{figure}

When we shave a moustache, we effectively map a point
$s \in \cell{G} \subseteq \Rspace^{nd}_\Delta$ to a point
$s' \in \cell{G'}$ in $(n-1)d$- or $(n-2)d$-dimensional configuration
space, depending on the case.
Doing this over all points $s$ and all moustaches of the corresponding
shortest networks,
we get multiple maps from higher- to lower-dimensional configuration spaces,
each skipping either $d$ or $2d$ dimensions.
We will exploit the fact that every point in an unambiguous cell
has plenty of preimages.
In other words, if $S'$ has an unambiguous shortest network,
then we can grow a moustache while preserving that the network is
shortest and unambiguous, of course now spanning a larger set $S$.
This property does not necessarily hold for sets $S'$ with
ambiguous shortest networks.
Indeed, the corners of a regular octagon have eight different shortest
networks \cite{DHW87},
and none of them permits the growth of a one- or two-sided
moustache such that the network remains shortest;
see Figure \ref{fig:octagon}.

\paragraph{Admissible extensions.}
Before proving the property for unambiguous shortest networks,
we introduce notation for adding a moustache.
Let $\Gnet{}' \colon G' \to \Rspace^d$ be a shortest network,
let $v$ be a degree-$1$ vertex with neighbor $w$ in $G'$,
and let $\varphi \in \Sspace^{d-1}$ be a direction.
To add a two-sided moustache with anchor $v$,
we require that the angle between $\varphi$ and the direction
$\omega = (\Gnet{}' (w) - \Gnet{}' (v)) / \dist{\Gnet{}' (w)}{\Gnet{}' (v)}$
be $120^{\circ}$.
Equivalently, we require $\scalprod{\varphi}{\omega} = - \frac{1}{2}$,
and if this equation is satisfied,
then we call $\varphi$ an \emph{allowed direction} for $v$.
For such a $\varphi$, there exists a unique
second direction $\psi \in \Sspace^{d-1}$ such that
$\scalprod{\psi}{\varphi} = \scalprod{\psi}{\omega} = - \frac{1}{2}$.
Choosing $r > 0$, we construct
$S = S' \setminus \{ \Gnet{}' (v) \} \cup \{ \Gnet{}' (v) + r \varphi ,
                                            \Gnet{}' (v) + r \psi \}$.
Similarly, we construct $G$ from $G'$ by adding vertices $p$ and $q$,
connecting both with new edges to $v$,
and we construct $\Gnet{} \colon G \to \Rspace^d$ such that the
restriction to $G'$ agrees with $\Gnet{}'$.
We call $r$ \emph{admissible} for $v$ and $\varphi$
if $\Gnet{}$ is the unique shortest network of $S$.

Adding a one-sided moustache is similar, except that $\varphi$
is an \emph{allowed direction}
if $\scalprod{\varphi}{\omega} \leq - \frac{1}{2}$,
and there is no second direction $\psi$ to worry about.

\begin{lemma}\label{lemma:admissible-extensions}
  {\sc (Existence of admissible extensions)}
  For each unambiguous shortest network $\Gnet{}' \colon G' \to \Rspace^d$,
  each degree-$1$ vertex $v$ in $G'$,
  and each allowed direction $\varphi \in \Sspace^{d-1}$,
  there exists $R > 0$ such that every $0 < r \leq R$
  is admissible for $v$ and $\varphi$.
\end{lemma}
\proof
 We consider the two-sided case and omit the proof of the easier, one-sided case.
 Assume without loss of generality that $G'$ is the type of $\Gnet{}'$,
 write $S'$ for the ordered set of points spanned by $\Gnet{}'$,
 and let $x = \Gnet{}' (v)$ be the image of the anchor.
 We prove that the set of admissible extensions is not empty.
 Given a positive admissible extension, Lemma \ref{lemma:trimming}
 then implies that all smaller extensions are also admissible,
 which gives the claimed statement.
 To get a contradiction, we assume that there is a degree-$1$ vertex
 $v$ of $G'$ and an allowed direction $\varphi$,
 such that no $r > 0$ is admissible for $v$ and $\varphi$.
 The plan is to construct two different locally minimal networks that
 are parametrized by the same tree.
 By Proposition \ref{prop:minimal-networks}, such networks do not exist,
 which will furnish the desired contradiction.

 The first locally minimal network is $\Gnet{} \colon G \to \Rspace^d$
 obtained by growing a two-sided moustache of sufficiently small
 size $R = r_1 > 0$ to avoid self-intersections.
 Assuming $G$ is its type, we get $G$ from $G'$ by adding vertices
 $p$ and $q$ and edges that connect them to $v$.
 To construct the second locally minimal network, we use the assumption
 that there is no positive admissible extension.
 Choosing an infinite sequence $r_i$ that vanishes in the limit,
 we set $S_i = S' \setminus \{x\} \cup \{x + r_i \varphi, x + r_i \psi\}$
 and write $\Hnet{i} \colon H_i \to \Rspace^d$ for the shortest
 network of $S_i$, letting $H_i$ be its type.
 Since we can select a converging subsequence,
 we may assume that $H = H_i$ for all $i$,
 and that $\Hnet{i}$ converges to a network $\Hnet{} \colon H \to \Rspace^d$.
 Since both $\Hnet{}$ and $\Gnet{}'$ are shortest, they are the same,
 only with different parametrizations.
 To analyze this difference, we consider the preimage,
 $H_x = \Hnet{}^{-1} (x)$.
 \begin{figure}[hbt]
  \centering \resizebox{!}{0.8in}{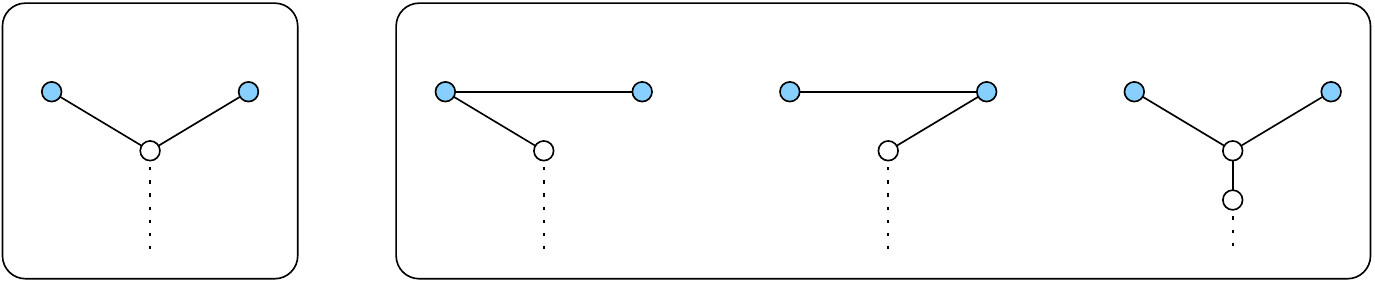}
  \caption{Possible subtrees $H_x$.
    The leftmost is a subtree of $G$, while the other three
    exhaust the possibilities for $H$.
    The rightmost subtree can be used to reparametrize the others.}
  \label{fig:preimages}
 \end{figure}
 We note that $p$ and $q$ are
 the only terminal vertices whose images converge to $x$.
 All interior vertices of $H$ have degree $3$, but one of them can have
 degree less than $3$ in $H_x$, namely the endpoint of the
 edge that connects $H_x$ to the rest of the tree.
 There cannot be a second such vertex for this would contradict that
 $\Hnet{}$ has a contracting parametrization
 as implied by Lemma \ref{lemma:convergence}.
 Since $H_x$ has at most three vertices of degree
 less than $3$, we are left with the choices shown in Figure \ref{fig:preimages}.
 The tree on the very right can be used to reparametrize the other three.
 We can therefore reparametrize
 $\Gnet{}$ and $\Hnet{1}$ with the same tree.
 Both networks are locally minimal and span the same set, $S_1$,
 which contradicts that they are different.
\eop

\paragraph{Lower semi-continuity.}
Fixing $G'$, we may consider the shortest networks of all
$s' \in \ucell{G'}$ and grow a moustache from the fixed anchor, $v$,
in each.
Within each such network, we have a direction
$\omega (s')$ determined by the image of $v$ and the image of
its neighbor, $w$.
By Lemma \ref{lemma:convergence}, $\omega$ depends continuously on $s'$.
We can therefore choose an allowed direction, $\varphi$,
that depends continuously on $s'$.
The second direction, $\psi$, needed to grow a
moustache is uniquely determined by $\omega$ and $\varphi$.
Let $R_1 \colon \ucell{G'} \to \Rspace$
be defined by mapping $s'$ to the minimum of $1$ and
the supremum of the thresholds $R$ for which
Lemma \ref{lemma:admissible-extensions} holds.
We do not know whether this function is continuous,
but we can prove lower semi-continuity, which is sufficient
to prove that the infimum of $R_1$
over a compact set of points $s'$ is positive.
\begin{lemma}\label{lemma:lower-semi-continuity}
  {\sc (Lower semi-continuity)}
  Let $s' = \lim_{i \to \infty} s_i'$, all points in a common
  unambiguous cell.
  Then $\liminf_{i \to \infty} R_1 (s_i') \geq R_1 (s')$.
\end{lemma}
\proof
 We consider the two-sided case and omit the proof of the easier, one-sided case.
 To get a contradiction, we suppose that
 $\liminf_{i \to \infty} R_1 (s_i') \leq R_1 (s') - 2 \delta$
 for some $\delta > 0$.
 Choose a converging sequence $r_i$ such that
 $R_1 (s_i') < r_i < R_1 (s') - \delta$ for each $i$.
 Writing $S_i'$ for the set of points that corresponds to $s_i'$
 and $x_i$ for the image of the anchor,
 we introduce
 $S_i = S_i' \setminus \{x_i\} \cup \{ x_i + r_i \varphi , x_i + r_i \psi \}$
 and let $\Gnet{i} \colon G \to \Rspace^d$ be the network
 obtained by adding the moustache to the shortest network of $S_i'$.
 Writing $G'$ for the type of that shortest network,
 we get the type $G$ of $\Gnet{i}$ from $G'$
 by adding vertices $p$ and $q$ and connecting them with edges to the anchor.
 Because $R_1 (s_i') < r_i$, the constructed network
 $\Gnet{i}$ is either not a shortest network or it is not unique.
 In either case, there exists a different shortest network
 $\Hnet{i}$ spanning the same ordered set $S_i \subseteq \Rspace^d$.
 Since we can choose subsequences, we may assume that the
 $\Hnet{i}$ all have the same type $H$.
 By Lemma \ref{lemma:convergence}, $\Hnet{i}$ converges to a shortest
 network $\Hnet{} \colon H \to \Rspace^d$ spanning
 $S = \lim_{i \to \infty} S_i$.
 By construction, $0 < \lim_{i \to \infty} r_i < R_1 (s')$,
 so $S$ has a unique shortest network $\Gnet{} \colon G \to \Rspace^d$.
 This implies that $\Gnet{}$ and $\Hnet{}$ are the same.
 Networks of type $G$ can therefore be parametrized by $H$.
 But then $\Gnet{i}$ and $\Hnet{i}$ are different locally minimal
 networks parametrized by the same tree,
 which contradicts Proposition \ref{prop:minimal-networks}.
\eop

\section{Proofs}
\label{sec4}

With the preparations in Section \ref{sec3}, we are now ready
to prove Theorems \ref{thm:connected-ucells} and \ref{thm:connected-cells},
which have been formally stated at the end of Section \ref{sec2}.
In addition, we illustrate the limitation of connectivity with
an example.

\paragraph{Proof of Theorem \ref{thm:connected-ucells}.}
Recall that this theorem asserts that unambiguous cells in the
decomposition of the configuration space of $n$ points
in $\Rspace^d$ are path-connected for all $d \geq 2$.
We use induction over the number of points.
For $n = 1$, there is only the trivial Steiner minimal tree,
whose unambiguous cell is the entire $\Rspace^d$, which is path-connected.
For $n = 2$, there is again only one Steiner minimal tree
consisting of two points connected by a line segment.
The unambiguous cell is the entire $\Rspace^{2d}$ minus the diagonal,
which is a $d$-dimensional plane.
This space has the homotopy type of the $(d-1)$-dimensional sphere.
It is path-connected for all $d \geq 2$.
This establishes the claim for $n = 1,2$ and $d \geq 2$.

For the inductive step, we consider two points, $s_0$ and $s_1$,
in an unambiguous cell $\ucell{G} \subseteq \Rspace^{(n+1)d}_\Delta$,
and prove that there exists a path $s \colon [0,1] \to \ucell{G}$
with $s(0) = s_0$ and $s(1) = s_1$.
Choose a moustache in $G$, shave it, and write $G'$ for the resulting
smaller tree.
This construction establishes a map from $\ucell{G}$
to $\ucell{G'}$, as described in Section \ref{sec3}.
Note that $\ucell{G'}$ is contained in $\Rspace^{nd}_\Delta$
or in $\Rspace^{(n-1)d}_\Delta$, depending on the choice of moustache.
By induction, $\ucell{G'}$ is path-connected, which implies
the existence of a path $s' \colon [0,1] \to \ucell{G'}$
such that $s'(0)$ and $s'(1)$ are images of $s_0$ and $s_1$
under the map.
We are going to lift the path to $\ucell{G}$ taking inverses.
Lemma \ref{lemma:admissible-extensions} states that we can
find a preimage of $s'(t)$ for every $t \in [0,1]$.
It remains to show that the preimages can be chosen in such a way
that they give a path $s \colon [0,1] \to \ucell{G}$.
But all the work needed to achieve this has already been done.
As explained before Lemma \ref{lemma:lower-semi-continuity},
we can choose directions $\varphi$ and $\psi$ that depend
continuously on $s'$.
It is easy to modify this direction field so that when we move
along $s'$, we connect the needed directions at the two ends.
Next, we set $r = \inf_{0 \leq t \leq 1} R_1 (s' (t))$
and note that $r > 0$ because $[0,1]$ is compact,
and $R_1$ is everywhere positive as well as lower semi-continuous;
see Lemma \ref{lemma:lower-semi-continuity}.
Moving along $s'$, we thus get a path $s_r \colon [0,1] \to \ucell{G}$
consisting of networks with moustaches of uniform size.
To get $s \colon [0,1] \to \ucell{G}$,
we first trim the moustache to go from $s_0$ to $s_r (0)$,
we second move from $s_r (0)$ to $s_r (1)$,
and we third enlarge the moustache to go from $s_r (1)$ to $s_1$.

\paragraph{Proof of Theorem \ref{thm:connected-cells}.}
Our second result asserts path-connectivity for cells, but only in
dimension $d = 2$ and for a full network in which all
terminal vertices have degree $1$.
Suppose $\Gnet{} \colon G \to \Rspace^2$ is a shortest such network
for a set $S$, let $G = (V \sqcup W, E)$ be its type,
and let $s_0 \in \cell{G} \setminus \ucell{G}$.
Since we know that $\ucell{G}$ is path-connected,
it suffices to construct a path $s \colon [0,1] \to \cell{G}$
that starts at $s(0) = s_0$ and ends at a point $s(1) \in \ucell{G}$.
We prove the existence of such a path by considering
the derivative of the length function.
Write $\mmm (u) \in \Rspace^2$ for the \emph{motion vector}
of $u \in V \sqcup W$, and define the \emph{direction vector}
by taking the sum of unit vectors in the directions of the
incident edges:
\begin{align}
  \ddd{\Gnet{}} (u)  &=  \sum_{\{u,v\} \in E}
    \frac{\Gnet{}(v) - \Gnet{}(u)}
         {\dist{\Gnet{}(v)}{\Gnet{}(u)}} .
\end{align}
The first derivative of the length function in the direction of the motion
is the sum of scalar products between motion and direction vectors.
Since the direction vector of every interior vertex vanishes,
this gives
\begin{align}
  \frac{\partial \Length{\Gnet{}}}{\partial t}  &=
    \sum_{u \in V} \scalprod{\mmm (u)}{\ddd{\Gnet{}} (u)} ;
\end{align}
see \cite[page 64]{IvTu94}.
We are interested in the special case in which
$\mmm (u) = \ddd{\Gnet{}} (u)$ for every $u \in V$.
Then $\partial \Length{\Gnet{}} / \partial t = n$,
the number of terminal vertices.
Returning to the task at hand, we consider a second shortest network,
$\Hnet{} \colon H \to \Rspace^2$ of $S$,
and we compare the two derivatives.
The vectors $\ddd{\Hnet{}} (u)$ have at most unit length.
This is immediate if the degree of $u$ in $\Hnet{}$ is $1$,
and it follows from the angle condition if the degree is $2$.
Hence, $\partial \Length{\Gnet{}} / \partial t
      - \partial \Length{\Hnet{}} / \partial t \geq 0$,
with equality iff $\ddd{\Gnet{}} (u) = \ddd{\Hnet{}} (u)$ for all terminal vertices.
If the difference is positive, then we move the points
$\Gnet{} (u)$ to $\Gnet{} (u) + r \ddd{\Gnet{}} (u)$
for $r > 0$ smaller than the length of the shortest edge of $\Gnet{}$.
After this motion, the adjusted network $\Gnet{}$ is shorter than
the adjusted network $\Hnet{}$.

It remains to show that vanishing difference is not possible.
As mentioned, this is equivalent to $\ddd{\Gnet{}} (u) = \ddd{\Hnet{}} (u)$
for all $u \in V$.
If all terminal vertices have degree $1$, then this implies that
$\Gnet{}$ and $\Hnet{}$ are codirected, which according to
Proposition \ref{prop:codirected} is impossible for two shortest networks
that are not the same.
Otherwise, we grow a new edge from every vertex $u$ whose degree
in $\Hnet{}$ is $2$ in the direction $- \ddd{\Gnet{}} (u)$.
We grow these edges both in $\Gnet{}$ and in $\Hnet{}$.
The two trees may no longer be shortest, but for sufficiently short
extensions, they are both locally minimal.
We thus get a contradiction from the stronger version of
Proposition \ref{prop:codirected} which asserts that no two locally minimal
trees can be codirected.

\begin{figure}[hbt]
  \centering
  \includegraphics[scale=0.7]{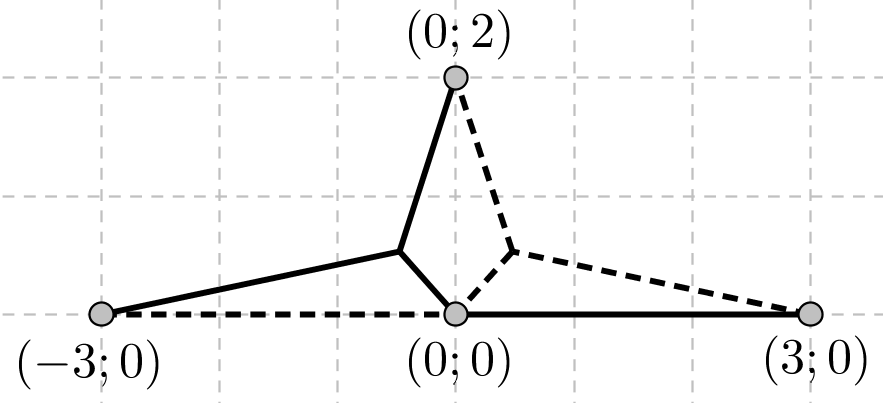} \hfill
  \includegraphics[scale=0.7]{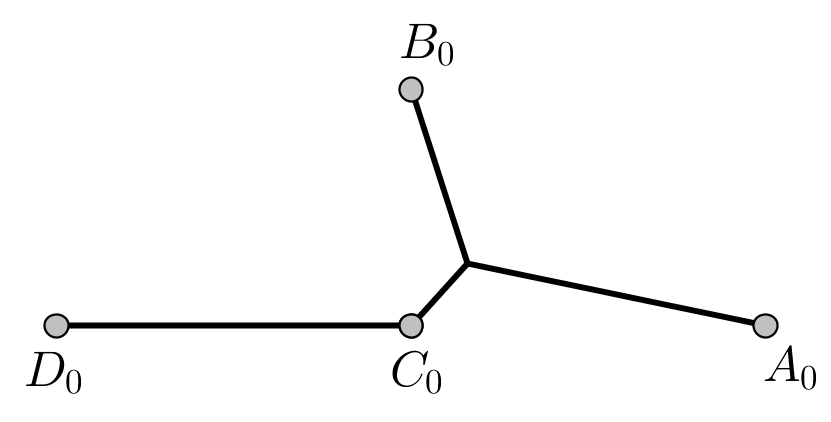} \hfill
  \includegraphics[scale=0.7]{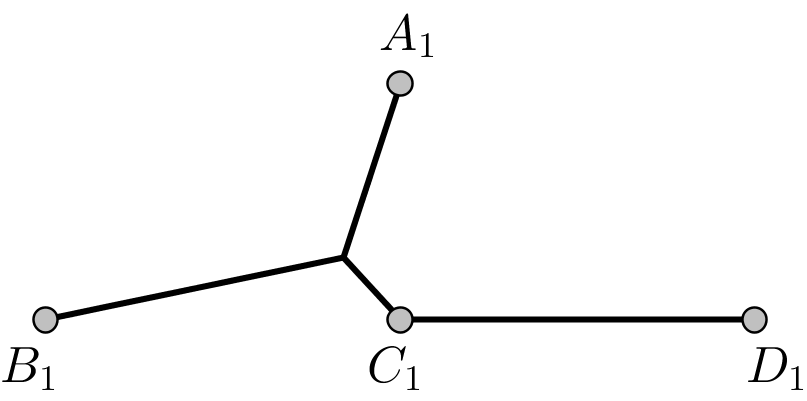}
  \caption{\emph{Left:} the two Steiner minimal trees of four points in the plane.
    \emph{Middle} and \emph{right:}  the two ambiguous shortest networks
    that cannot be connected by a path of like networks.}
  \label{fig:ambiguous}
\end{figure}

\paragraph{Example.}
To complement our theorems, we show that ambiguous subsets of
cells are not necessarily connected.
Specifically, we exhibit an ambiguous shortest network of
type $G$ such that $\cell{G} \setminus \ucell{G}$ is not connected.
The network spans the ordered set of $n = 4$ points in $\Rspace^2$
drawn in Figure \ref{fig:ambiguous}.
The two trees superimposed on the left are the only two Steiner minimal trees
of the four points.
Indeed, there cannot be a full Steiner minimal tree with two interior vertices
because the four points do not form a strictly convex quadrilateral;
see \cite[Lemma 5]{Pol78},
and the remaining types with one or zero interior
vertices give longer networks.

For the next step, we choose two orderings of the four points such that
the first Steiner minimal tree for the first ordering has the same type
as the second Steiner minimal tree for the second ordering;
see Figure \ref{fig:ambiguous} in the middle and on the right.
We claim that there is no path from the first to the second tree along
which the Steiner minimal tree remains of the same type and ambiguous at all times.
To derive a contradiction, we suppose there is such a path
$s \colon [0,1] \to \cell{G} \setminus \ucell{G}$.
Consider the counterclockwise angle, $\alpha_t$,
from the edge connecting $C_t$ to $D_t$ to the second edge incident on $C_t$.
Since we have a Steiner minimal tree at all times,
we have $\frac{2 \pi}{3} \leq \alpha_t \leq \frac{4 \pi}{3}$ for all $t$.
Taking a look at Figure \ref{fig:ambiguous}, we see that
$\alpha_0 > \pi$ and $\alpha_1 = 2 \pi - \alpha_0 < \pi$.
Since $\alpha_t$ depends continuously on $t$, there is a value $0 < \tau < 1$
such that $\alpha_{\tau} = \pi$.
At $t = \tau$, the Steiner minimal tree is determined by the points
$A_\tau, B_\tau, D_\tau$.
However, the Steiner minimal tree of three points cannot be ambiguous.
This contradicts that $S$ is a path inside $\cell{G} \setminus \ucell{G}$.
We note that this example does not extend to $d \geq 3$ dimensions.

\section{Discussion}
\label{sec5}

This paper continues the study of the decomposition of the
configuration space of $n$ points in $\Rspace^d$ as defined
by the combinatorial types of the Steiner minimal trees,
which was initiated in \cite{IvTu06}.
We consider the ordered setting and our main result is a proof
that the cells consisting of configurations with unambiguous
Steiner minimal trees are path-connected.
There are many questions that remain open.
\begin{itemize}
  \item  Can our partial results for cells be extended to all
    configurations in $\Rspace^2$ and to higher dimensions?
  \item  Is it true that the set of points in $\Rspace^{nd}_\Delta$
    with ambiguous shortest networks has measure zero?
  \item  Is the non-empty intersection of two cells in $\Rspace^{nd}_\Delta$
    necessarily path-connected?
  \item  Is the union of sets $\cell{G} \setminus \ucell{G}$
    in $\Rspace^{nd}_\Delta$ over all types $G$ path-connected?
\end{itemize}
Besides the ordered setting considered in this paper, it would be interesting
to address the same questions for unordered point sets.
The additional symmetries complicate the global topology \cite{Bir74},
and more sophisticated notions of connectivity as offered by homology
and homotopy groups seem appropriate.

\paragraph{Acknowledgments.}
{\small The authors thank A.\,O.~Ivanov, Z.\,N.~Ovsyannikov,
and A.\,A.~Tuzhilin for insightful discussions on the contents of this paper.}


\newpage \appendix
\section{Notation}

\begin{table}[hbt]
  \centering
  \begin{tabular}{ll}
    $S \subseteq \Rspace^d$
      &  ordered set of $n$ points in Euclidean space            \\
    $G = (V \sqcup W, E)$
      &  partially ordered connected graph                       \\
    $\Gnet{} \colon G \to \Rspace^d, G_0$
      &  network, graph; type                                    \\
    $\Hnet{} \colon H \to \Rspace^d, H_0$
      &  network, graph; type                                    \\
    $g, g_0, h$
      &  homeomorphisms                                          \\
                                                                 \\
    $s \in \Rspace^{nd}_\Delta$
      &  configuration, configuration space                      \\
    $\ucell{G} \subseteq \cell{G}$
      &  unambiguous cell, cell                                  \\
    $s \colon [0,1] \to \ucell{G}$
      &  path                                                    \\
    $t, \tau \in [0,1]$
      &  time parameter, moment                                  \\
                                                                 \\
    $\Gnet{i}, i \in \Nspace$
      &  infinite sequence, index                                \\
    $\Gnet{i} \to \Gnet{}, \Gnet{} = \lim \Gnet{i}$
      &  convergence, limit                                      \\
    $C; D \subseteq \Rspace^d$
      &  constant; bounded convex region                         \\
    $u, v, w; p, q$
      &  vertices                                                \\
    $\varphi, \psi, \omega \in \Sspace^{d-1}$
      &  directions                                              \\
    $r_i, r, R$
      &  positive extensions                                     \\
    $R_1 \colon \ucell{H} \to \Rspace$
      &  supremum extension function                             \\
                                                                 \\
    $\mmm (u), \ddd{\Gnet{}} (u)$
      &  motion vector, direction vector                         \\
    $\Length{\Gnet{}}, \partial \Length{\Gnet{}} / \partial t$
      &  length, derivative                                      \\
    $\alpha, \alpha_t$
      &  angles
  \end{tabular}
  \caption{Notation for geometric concepts, sets, functions,
    vectors, variables used in this paper.}
  \label{tbl:Notation}
\end{table}

\section{Dependence}
\begin{figure}[hbt]
 \centering \resizebox{!}{2.0in}{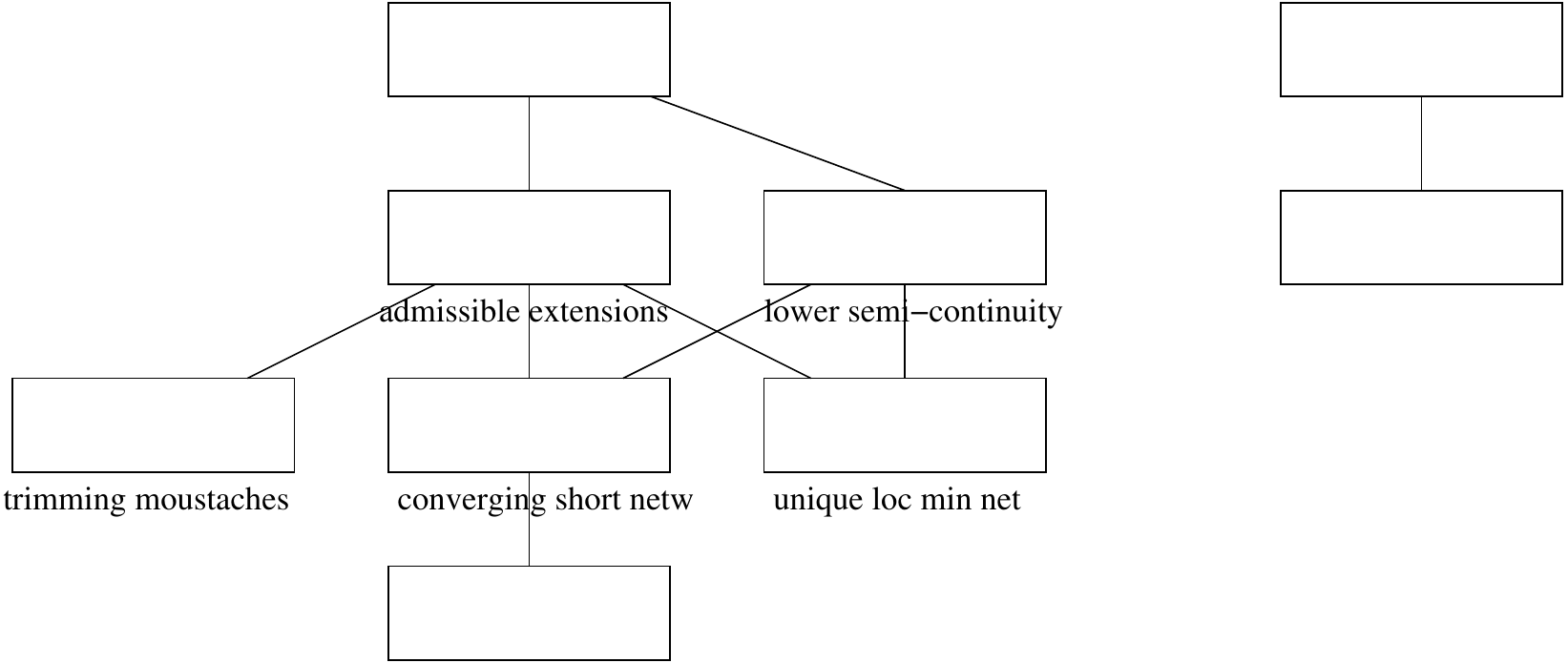}
 \caption{Dependence of Lemmas and Theorems from top to bottom.}
 \label{fig:dependence}
\end{figure}

\end{document}